\documentclass[letterpaper, 10 pt, conference]{ieeeconf}  %

\IEEEoverridecommandlockouts                              %

\usepackage{tikz-cd}
\tikzcdset{/tikz/commutative diagrams/background color=red}
\tikzcdset{row sep/normal=1cm}

\usepackage{graphicx} %
\usepackage{epsfig} %
\usepackage{amsmath} %
\usepackage{amssymb}  %
\usepackage{booktabs}
\usepackage{makecell}
\usepackage{hyperref}
\usepackage{subcaption}
\usepackage{optidef}
\usepackage{comment}

\newcommand{\R}{\mathbb{R}}

\newcommand{\ind}[1]{{_{\mathrm{#1}}}}
\newcommand{\diag}[1]{\text{\normalfont diag}\left( {#1} \right)}

\newcommand{\acados}{\textup{\texttt{acados}}}

\newcommand{\Set}[2]{\{\, #1 \mid #2 \,\}}

\DeclarePairedDelimiter\abs{\lvert}{\rvert}%
\DeclarePairedDelimiter\norm{\lVert}{\rVert}%

\makeatletter
\let\oldabs\abs
\def\abs{\@ifstar{\oldabs}{\oldabs*}}
\let\oldnorm\norm
\def\norm{\@ifstar{\oldnorm}{\oldnorm*}}
\makeatother

\DeclareMathSymbol{\sm}{\mathbin}{AMSa}{"39}

\newcommand{\eye}{%
        \text{\usefont{U}{bbold}{m}{n}1}%
}
\newcommand{\eyed}[1]{%
        \eye_{#1}%
}

\newcommand{\tran}{^\top}

\newcommand{\nx}{{n_{\mathrm{x}}}}

\DeclareMathSymbol{\shortminus}{\mathbin}{AMSa}{"39}
\newcommand{\statecov}{P}

\newcommand{\dtotal}[2]{\frac{\mathrm{d} #1}{\mathrm{d} #2}}
\newcommand{\dpartial}[2]{\frac{\partial #1}{\partial #2}}
\newcommand{\zoro}{zoRO}
\newcommand{\oldzoro}{\mbox{zoRO-21}}
\newcommand{\newzoro}{\mbox{zoRO-24}}

\newcommand{\todofinal}[1]{}

\title{\LARGE \bf
Efficient Zero-Order Robust Optimization for Real-Time Model Predictive Control with \acados
}

\author{Jonathan Frey$^{1,2}$, Yunfan Gao$^{1,3}$, Florian Messerer$^{1}$, Amon Lahr$^4$, Melanie Zeilinger$^{4}$,
Moritz Diehl$^{1,2}$%
\thanks{$^{1}$Department of Microsystems Engineering (IMTEK), University Freiburg, 79110 Freiburg, Germany
{\tt\small \{name.surname\}@imtek.uni-freiburg.de}
}
\thanks{$^{2}$Department of Mathematics, University Freiburg, Germany
}%
\thanks{$^{3}$Robert Bosch GmbH, Corporate Research, Stuttgart, Germany
{\tt\small yunfan.gao@de.bosch.com}
}
\thanks{$^4$Institute for Dynamic Systems and Control, ETH Zurich, 8092 Zurich, Switzerland}
\thanks{This research was supported by DFG via Research Unit FOR 2401, project 424107692 on Robust MPC and 525018088, by BMWK via 03EI4057A and 03EN3054B, and by the EU via ELO-X 953348.}
}

\begin{document}

\maketitle
\thispagestyle{empty}
\pagestyle{empty}

\begin{abstract}
Robust and stochastic optimal control problem (OCP) formulations allow a systematic treatment of uncertainty, but are typically associated with a high computational cost.
The recently proposed zero-order robust optimization (zoRO) algorithm mitigates the computational cost of uncertainty-aware MPC by propagating the uncertainties outside of the MPC problem.
This paper details the combination of zoRO with the real-time iteration (RTI) scheme and presents an efficient open-source implementation in \texttt{acados}, utilizing \texttt{BLASFEO} for the linear algebra operations.
In addition to the scaling advantages posed by the zoRO algorithm, the efficient implementation drastically reduces the computational overhead, and, combined with an RTI scheme, enables the use of tube-based MPC for a wider range of applications.
The flexibility, usability and effectiveness of the proposed implementation is demonstrated on two examples.
On the practical example of a differential drive robot, the proposed implementation results in a tenfold reduction of computation time with respect to the previously available zoRO implementation.

\end{abstract}
\section{Introduction}
A dedicated, explicit treatment of uncertainties allows practical model predictive control (MPC) applications to avoid heuristic and overly conservative safety margins and to harness unused optimality potential.
A wide variety of uncertainty-aware OCP formulations exists, such as min-max MPC, scenario-tree MPC and tube-based MPC with tubes of various shapes~\cite{Mayne2011},~\cite{Kouzoupis2019},~\cite{Telen2015},~\cite{Gillis2013}.
This paper focuses on problem formulations where the uncertainties are represented by ellipsoidal tubes in the robust case, and independent distributed noise in the stochastic case.
In the case of linear dynamics and constraints, constraint tightenings can be precomputed and such problems can be solved at no additional computational cost compared to their nominal correspondents.
In the case of nonlinear dynamics, it is sometimes possible to design constraint tightenings offline by scaling the uncertainty set~\cite{Koehler2018}.
However, the state-independent over-approximation of the uncertainty may lead to significant conservatism of the controller, especially when the shape of the actual uncertainty set strongly varies along the prediction horizon.

In this work, we thus consider constraint tightenings based on an approximate uncertainty propagation performed online. Treating the nonlinear form of such problems conventionally with standard OCP solvers requires a state augmentation that is quadratic in the original state dimension.
Since the computational cost of most OCP solvers is cubic in the state dimension, this complicates the exact treatment of tube-based OCPs in real-time MPC applications with an overall computational burden growing with the sixth power of the state dimension.
As a remedy, the recently proposed zero-order robust optimization (zoRO) algorithm allows one to treat these OCPs at a computational complexity of the order corresponding to the one of a nominal OCP~\cite{Zanelli2021, Feng2020}, while returning a suboptimal, yet feasible, point at convergence~\cite{Zanelli2021}.
This algorithmic novelty has brought the application of tube-based MPC with online uncertainty propagation on real systems within reach, as zoRO has been applied successfully on the experimental setup of a differential drive robot~\cite{Gao2023}.
Furthermore, the zoRO algorithm has been extended in~\cite{Lahr2023} to incorporate learning-based model uncertainty, e.g., Gaussian process-based MPC~\cite{Hewing2020} or Bayesian last-layer networks~\cite{Fiedler2022}.
In addition, the problem formulation allows one to incorporate precomputed linear feedback laws.

Existing zoRO implementations have been written in Python and are limited to their specific use-cases~\cite{Zanelli2021, Gao2023, Lahr2023}.
In contrast to this, we present an efficient, flexible, open-source \texttt{C} implementation of zoRO, which can be conveniently used from Python.
The high-performance implementation utilizes intermediate results from the~\acados~SQP solver and performs the uncertainty propagation as well as backoff computation using \texttt{BLASFEO} for the linear algebra operations~\cite{Frison2018}, leading to significant computational speedups compared to the previous implementations.

In the following, we give a brief introduction into the OCP formulation (Section~\ref{sec:robust_stochastic_oc}) and the zoRO algorithm (Section~\ref{sec:zoro}), before presenting the main features of the high-performance implementation (Section~\ref{sec:acados_zoro}).
Finally, Section~\ref{sec:experiments} presents numerical experiments and Section~\ref{sec:conclusion} concludes the paper.

\newpage
\section{Robust \& Stochastic Optimal Control}
\label{sec:robust_stochastic_oc}

\vspace{-.0cm}
Let us regard optimal control problems~(OCPs) of the form
\begin{mini!}|s|%
	{\substack{x_0, \dots, x_N, \\ u_0, \dots, u_{N \shortminus 1}\\  \statecov_{0}, \dots, \statecov_{N}}}
	{\sum_{k=0}^{N-1} l(u_k, x_k) + M(x_N)}
	{\label{eq:robustMPC}}{}
	\addConstraint{x_0}{=\bar{x}_{0}}{}
	\addConstraint{\statecov_0}{= \bar{\statecov}_{0}}{}
	\addConstraint{x_{k+1}}{=\psi_k(x_k, u_k, 0)\label{eq:nominal_dyn}}{}
	\addConstraint{\statecov_{k+1}}{=\Phi_k(x_k, u_k, \statecov_k)\label{eq:uncertainty_dyn}}{}
	\addConstraint{0}{\geq h_k(x_k, u_k) + \beta_k(x_k, u_k, \statecov_k \label{eq:ineq_path})}{}
	\addConstraint{0}{\geq h_N(x_N) +   \beta_N(x_N, \statecov_N \label{eq:ineq_terminal})}
	\addConstraint{}{\quad\quad \mathrm{with} \quad k = 0,\dots, N-1 \nonumber}%
\end{mini!}
where $x_k \in \R^{n_x}$ are the state variables for $k=0, \dots, N$, and $u_k \in \R^{n_u}$ are the controls for $k=0, \dots, N-1$.
The nonlinear, discrete-time system dynamics are given by $\psi_k: \R^{n_x}\times \R^{n_u} \times \R^{n_w} \to \R^{n_x}$ and additionally depend on an uncertain variable $w_k \in \R^{n_w}$.
The resulting uncertain trajectory is represented by the nominal trajectory $ (x_0, u_0, \dots, u_{N-1}, x_N) $ and the matrices $\statecov_0, \dots, \statecov_N \in \mathbb{R}^{n_x \times n_x}$.
A simulation with zero noise characterizes the nominal trajectory~\eqref{eq:nominal_dyn}; equation~\eqref{eq:uncertainty_dyn} describes the uncertainty dynamics.
Here, the positive (semi)definite matrices $\statecov_k$ represent uncertainty in state space and have the meaning of covariance matrices or ellipsoidal shape matrices.

This formulation covers uncertainty-aware NMPC for both a stochastic and robust setting.
In the stochastic setting, each $w_k$ follows an independent normal distribution with zero mean and covariance $W_k$, e.g. $w_k\sim \mathcal{N}(0, W_k)$.
Here, $x_k$ and $\statecov_k$ parameterize mean and variance of a stochastic state $\chi_k$, usually approximated as $ \chi_k \sim \mathcal{N}(x_k, \statecov_k)$, which approximates the exact state distribution over the horizon for fixed controls.

In the robust setting, the assumption is that the noise variables are contained within an ellipsoidal set
\begin{align}
	(w_0, \dots, w_{N-1}) \in \mathcal{E}\left(0, \mathrm{blkdiag}( W_0,\dots, W_{N-1}) \right).
\end{align}
Here, $\mathcal{E}(q, Q) := \Set{x\in \R^{n_q}}{(x-q)^{\top} Q^{-1} (x-q) \leq 1} $ denotes an ellipsoid with center $q \in \R^{n_q}$ and shape described by a positive definite $Q\in \R^{n_q \times n_q}$.
This results in an ellipsoidal tube such that $\chi_k \in \mathcal{E}(x_k, P_k)$.
\\
The nominal initial state is $\bar{x}_0$ and has uncertainty described by $\bar{P}_0$, e.g. the mean and covariance of the current state estimate in the stochastic setting.
The inequality constraints~\eqref{eq:ineq_path} and~\eqref{eq:ineq_terminal} consist of the nominal term $ h_k (\cdot)$ and an additional backoff $ \beta_k(\cdot) $ to account for the uncertainty,
representing chance constraints in the stochastic setting and worst-case constraints in the robust setting~\cite{Ben-Tal2009}.

The uncertainty propagation can be written as
\begin{align}
	\Phi_k(P_k&, x_k, u_k) \label{eq:uncertainty_prop}\\
	&= (A_k + B_k K_k) P_k (A_k + B_k K_k)^\top + G_k W_k G_k^\top, \nonumber
\end{align}
where $ A_k := \dpartial{\psi_k}{x}(x_k, u_k, 0) $, and $ B_k := \dpartial{\psi_k}{u}(x_k, u_k, 0) $.
The matrices $G_k = \dpartial{\psi_k}{w}(x_k, u_k, 0) $
are constant in case of additive noise.
Additionally, it is possible to include precomputed linear feedback gains $K_k\in\R^{n_u \times n_x}$ to reduce conservatism~\cite{Mayne2011}.
\\
The backoff term for a constraint component $ h_{k, i}(\cdot)$ of $ h_k\!:\! \R^{n_x} \!\times \! \R^{n_u} \!\to \!\R^{n_{h,k}} $ for $ i \in \{1, \dots, n_{h,k}\} $ is given by
\begin{align}
\beta_{k,i}&(x_k, u_k, \statecov_k):= \label{eq:backoffs} \\
& \gamma
\sqrt{\nabla h_{k,i}(x_k, u_k)^{\top}\!
	\left[\begin{array}{c}
	\eyed{\nx}\!\\ K\!
	\end{array} \right]  \!P_k
	{\left[\begin{array}{c}
		\eyed{\nx}\\ K
		\end{array} \right]}^{\!\top}\nabla h_{k,i}(x_k, u_k)}. \nonumber
\end{align}
In the robust setting, the backoff factor $\gamma$ equals $1$, while in the stochastic setting, $\gamma$ is the number of standard deviations that the nominal trajectory should maintain a distance from the nominal bounds.
Depending on the interpretation of $x_k, \statecov_k$ as a normal distribution or the first two moments of a general distribution, the choice of $\gamma$ corresponds to different probabilities of constraint satisfaction, which can be derived via Chebyshev's inequality or the inverse normal cumulative density function, c.f.~\cite{Heirung2018}.

\section{Zero-Order Robust Optimization (zoRO)}
\label{sec:zoro}
The key idea of zoRO is to eliminate the uncertainty matrices from the OCP and solve reduced subproblems in the nominal variables only,
\begin{mini!}|s|
	{\substack{x_0, \dots, x_N, \\ u_0, \dots, u_{N \shortminus 1}}}
	{\sum_{k=0}^{N-1} l(u_k, x_k) + M(x_N)}
	{\label{eq:zoro_ocp}}{}
	\addConstraint{x_0}{=\bar{x}_{0}}{}
	\addConstraint{x_{k+1}}{=\psi_k(x_k, u_k, 0)}{}
	\addConstraint{0}{\geq h_k(x_k, u_k) + \hat{\beta}_k}{}
	\addConstraint{0}{\geq h_N(x_N) +  \hat{\beta}_N}
	\addConstraint{}{\quad\quad \mathrm{with} \quad k = 0,\dots, N-1 \nonumber,}%
\end{mini!}
where the the values $\hat{\beta}_k$ approximate the backoff terms $\beta_k(x_k, u_k, \statecov_k)$. For a detailed derivation, we refer the interested reader to~\cite{Zanelli2021, Feng2020}.

The zoRO algorithm alternates the following two steps:
\begin{enumerate}
	\item approximate the backoff terms using the current guess of the nominal trajectory, by performing the uncertainty propagation~\eqref{eq:uncertainty_prop} and the backoff computation~\eqref{eq:backoffs}
	\item solve subproblem~\eqref{eq:zoro_ocp} approximately
\end{enumerate}
\todofinal{Check the indents of each paragraph such that you have indents where you would like to have them (also in previous paragraphs)}
Note that the accuracy up to which the subproblems are solved is an implementation choice and depends on the solver used to tackle the subproblems~\eqref{eq:zoro_ocp}.
For example, the subproblems might be solved to convergence using an interior point method, or with a limited number of SQP iterations, c.f. Section~\ref{sec:zoro_SQP}.

\subsection{Convergence properties} \label{sec:convergence}

Comparing the solution of original problem~\eqref{eq:robustMPC} with a Newton-type optimization method and the alternation performed by zoRO, it can be seen that, while in the original problem the uncertainty matrices correspond to optimization variables, in the reduced subproblem they are inserted as fixed parameters.
Thus, it is possible to reformulate the zoRO algorithm as an inexact Newton-type method, which employs a tailored Jacobian approximation at each iteration, namely, by neglecting the sensitivities of the uncertainty matrices with respect to the state and input variables~\cite{Zanelli2021,Feng2020}.

By standard arguments for Newton-type optimization, see e.g.~\cite{Bock2007}, the zoRO algorithm hence returns a suboptimal, yet feasible, solution of the original problem.

Moreover, due to the structure of the tailored Jacobian approximation, the approximation error scales with the magnitude of the uncertainty, leading to favorable convergence properties for sufficiently small level of uncertainty, shown in~\cite{Zanelli2021}.
First, the incurred suboptimality approaches zero, i.e., the zoRO algorithm recovers the optimal solution to~\eqref{eq:robustMPC} in the limit for vanishing uncertainties.
Second, for sufficiently small uncertainties, zoRO converges linearly provided that the unmodified SQP algorithm converges.
While the results in~\cite{Zanelli2021} hold for the case of constant additive uncertainties~$W$, the second convergence result has been extended to the case of state- and input-dependent uncertainties~$W(y)$ in~\cite{Lahr2023}.

Note that it is possible to compensate for the suboptimality that the tailored Jacobian approximation incurs by means of an adjoint correction~\cite{Bock2007,Wirsching2006}, as done in~\cite{Feng2020}, at the additional expense of computing the neglected sensitivities in the direction of the corresponding Lagrange multipliers with the backwards-mode of automatic differentiation.
In this paper, we do not address this variant due to the higher computational cost.

\subsection{Remark on Nonlinear Constraints}

Previous works have used the term zoRO for two slightly different variants of the algorithm.
In~\cite{Gao2023}, the NLP with fixed backoff terms~\eqref{eq:zoro_ocp} is solved with an SQP algorithm.
On the other hand,~\cite{Zanelli2021, Feng2020, Lahr2023} directly take the SQP perspective and suggest solving quadratic subproblems with the constraints
\begin{align}
	\label{eq:sqp_backoffs}
0 &\geq h_{k,i}(\bar{y}_k) + \beta_{k,i}(\bar{y}_k, \bar{P}_k) + \nabla_y h_{k,i}(\bar{y}_k) \Delta y_k \\
& +\nabla_y \beta_{k,i}(\bar{y}_k, \bar{P}_k)^\top \Delta y_k +\nabla_P \beta_{k,i}(\bar{y}_k, \bar{P}_k)^\top \Delta P_k \nonumber
\end{align}
see~\cite[Eq. (14)]{Zanelli2021}, where $ \bar{y}_k = (\bar{x}_k, \bar{u}_k) $ and $\bar{P}_k$ is the current linearization point.

Thus, compared with~\cite{Zanelli2021, Feng2020, Lahr2023}, performing a single SQP iteration on~\eqref{eq:zoro_ocp} corresponds to additionally neglecting the derivative of the backoff with respect to the nominal variables $\nabla_y \beta_{k,i}$.
Further, in~\eqref{eq:zoro_ocp}, the backoff $\beta_{k,i}$ is evaluated at $\bar{P}_k^+$, instead of linearizing $h_{k,i}$ with respect to $P_k$, as in~\eqref{eq:sqp_backoffs}.

To compute the neglected derivatives, on the one hand, the term $\nabla_P \beta_{k,i} \Delta P_k$ has to be precomputed and added, respectively subtracted, on top of the back-off term, since $P_k$ is merely a parameter of the reduced subproblems.
Thereby, the difference $\Delta P_k$ is computed as $\Delta P_k = \bar{P}_k^+ - \bar{P}_k$, where $\bar{P}^+_k$ denotes the uncertainty matrix after the performing the uncertainty propagation at~$\bar{y}_k$.
On the other hand, for the term $\nabla_y \beta_{k,i}(\bar{y}_k, \bar{P}_k)^\top \Delta y_k$ the derivative of the backoff terms $\beta_{k,i}$ with respect to the nominal variables $x_k, u_k$ are needed.
Those are
\begin{align}
	\label{eq:backoff_nonlinear_constraints}
	\nabla_{y} \beta_{k,i} &=  \frac{\gamma^2}{\beta_{k,i}} \cdot \\
	&\left( \nabla_y h_{k,i}(\cdot, \cdot)^\top \!
	\left[\begin{array}{c}
		\eyed{\nx}\\ K
	\end{array} \right]  \!P_k
		{\left[\begin{array}{c}
			\eyed{\nx}\\ K
			\end{array} \right]}^{\!\top\!} \!
			\nabla_{yy} h_{k,i}(\cdot, \cdot)
			\right)^{\!\top}\!, \nonumber
\end{align}
which vanish for linear constraints $h_{k,i}$.
Neglecting those terms results again in an inexact Jacobian.

We note that the convergence results of~\cite{Zanelli2021} for the zoRO algorithm still hold under the same assumptions when neglecting these additional derivatives.
This can be seen by the fact that the additionally neglected parts of the Jacobian also scale with the magnitude of the uncertainty.
Thus, they similarly recover the unmodified Jacobian in the limit for vanishing uncertainties and, ultimately, the convergence properties outlined in Section~\ref{sec:convergence}.

For computational efficiency and ease of implementation, the remainder of this work considers the zoRO version that neglects the terms in~\eqref{eq:backoff_nonlinear_constraints}.

\section{\acados~zoRO SQP implementation}
\label{sec:acados_zoro}
In this section, we detail the ingredients for an efficient zoRO implementation with SQP and RTI in \texttt{acados}~\cite{Verschueren2021}.

\subsection{Real-time iterations}
The real-time iteration (RTI) scheme~\cite{Diehl2005c} allows one to split an SQP iteration into a preparation and a feedback phase.\\
The preparation phase is carried out based on the initial guess and evaluates the nonlinear functions and its (approximate) derivatives to form the (approximate) Hessian of the Lagrangian.
The Hessian can then be used to prepare the linear algebra operations for the feedback phase.
However, in case of an interior-point QP solver, like \texttt{HPIPM}, the scope of operations during the preparation phase is limited.

The feedback phase takes the best available estimate of the current state value $\bar{x}_0$, evaluates the remaining parts of the Lagrange gradient, and solves the QP subproblem.

We note that typically, only the component corresponding to the constraint $\bar{x}_0 - x_0$ is evaluated in this phase~\cite{Diehl2005c}.
However, to allow updating all constraint bounds, such as the backoffs in the zoRO algorithm, the \texttt{acados} RTI implementation was adapted as follows:
It makes all constraint evaluations available in the preparation phase, and just subtract its bounds in the feedback phase, to complete the computation of the Lagrange gradient.
The only additional computational cost for this is performing additions corresponding to the number of constraints in the feedback instead of the preparation phase.

\subsection{zoRO with SQP \& RTI}
\label{sec:zoro_SQP}
When using an SQP solver, as is the focus in this paper, a common choice is to update the backoffs after each SQP iteration.
For an efficient implementation of zoRO-SQP, the current linearization of the dynamics and constraints, namely $A_k, B_k$ and $\nabla h_{k,i}(x_k, u_k) $ can be taken directly from the memory of the solver.
The backoffs can then be inserted into the subproblem by adapting the bounds of the constraints.

More specifically, when using the RTI variant, one can compute the backoff terms as part of the preparation phase, i.e. first perform the preparation step of~\eqref{eq:zoro_ocp} then update the backoff terms.
Once the new state estimate is available, use $\bar{x}_0$ to perform the feedback phase on the QP approximating the subproblem~\eqref{eq:zoro_ocp} and deploy the solution.
Note that, when employing the RTI variant of zoRO, the algorithm is equivalent to performing RTI and using the previous MPC solution to propagate the uncertainties and compute the backoffs.
Using the previous MPC solution for backoff computation with reduced computation time has also been performed in~\cite{Hewing2020}.

\subsection{Remark on non-tightened constraints}
In many practical OCP formulations, it is not desired to tighten all the constraints.
Thus, we define the index sets $\mathcal{I}_{k, \mathrm{tight}} \subseteq \{1,\dots, n_{h,k}\}$, which denote the subset of constraints that are tightened in~\eqref{eq:robustMPC}, i.e., for all other constraints $ h_{k,i}$ with $ i \notin \mathcal{I}_{k, \mathrm{tight}}$, the corresponding backoff is always zero.

Constraints that are not supposed to be tightened are for example bounds on slack variables.
Additionally, two inequalities can be used to encode an equality constraint, such as for the initial state constraint.
Furthermore, inequalities in \acados~are always formulated as two-sided, and if one of the bounds is a place holder, those do not need to be tightened.

\subsection{\acados~template interface}
An important aspect of the presented implementation is to specify information on the zoRO algorithm compactly in a predefined format, transfer it to \texttt{C} code and couple it with the \acados~solver.
This is achieved by utilizing the \acados~template interface workflow, which is briefly described in the following.

The \acados~high-level interfaces to Python and Matlab allow one to conveniently and compactly formulate an OCP and specify solver options.
The nonlinear problem functions can be formulated using \texttt{CasADi}~\cite{Andersson2019} symbolics which are generated as \texttt{C} code with the required derivatives using automatic differentiation.
The whole problem description is written to a \texttt{json}-file which is then used to render different templates.
Most importantly, a \texttt{C} file is generated which uses the \acados~shared library to create a solver specific to the problem formulation and loads the nonlinear functions.

\subsection{Generic Custom Updates}
In order to allow updating numerical data in the solver efficiently, i.e., without any calling functions outside the \texttt{C} stack, we added the option of using a custom \texttt{C} update function in between solver SQP solver calls.
This custom update function
\begin{itemize}
\item can access the solver;
\item can allocate its own memory and store a pointer to it in the solver \texttt{capsule}, implemented by calling an \texttt{\_init} function at the end of the creation process of the \acados~OCP solver and a \texttt{terminate} function which is responsible for freeing it;
\item is compiled together with the problem specific OCP solver;
\item can be conveniently used in the prototyping phase, as it is fully integrated in the Python interface and can simply be called using:\\
\texttt{AcadosOcpSolver.custom\_update()}.
\end{itemize}
\todofinal{Would be awesome if the intemize environment ends up in a single column in the end}
Moreover, the user is allowed to write information specific to their problem into the \texttt{json} file and use that to generate the custom update function from a custom template.
The \texttt{C} template implements the \zoro~update and allows the user to specify various options, as detailed in the next section.

We note that this feature can be useful to efficiently implement a variety of other advanced MPC schemes, other than zoRO, which require one to update parameters, constraint bounds, or other numerical data in the OCP solver in between SQP iterations.
Some examples include
\begin{itemize}
	\item the advanced-step RTI variants~\cite{Nurkanovic2019a},
	\item an efficient algorithm for robust MPC with optimal linear feedback~\cite{Messerer2021},
	\item custom backoff computations based on, e.g., poly- or zonotopic tubes,
	\item state-dependent modelling, e.g. via splines, without explicit incorporation (and therefore differentiation) in the nonlinear functions of the SQP scheme.
\end{itemize}

\subsection{Template-based \zoro~implementation}
The uncertainty propagation and the constraint tightening of the \zoro~algorithm is implemented in a template-based \texttt{C} function using \texttt{BLASFEO}~\cite{Frison2018} for all the linear algebra operations.
The constraint and model linearizations, i.e., $\nabla h_k, A_k, B_k$, are obtained from the \acados~solver, using the \texttt{C} interface, while linear constraints and simple bounds on states and controls are handled efficiently.

We added the Python class $\texttt{ZoroDescription}$ to conveniently specify the zoRO specific information. Specifically, this contains:
\begin{itemize}
	\item the initial uncertainty matrix $\bar{P}_0$,
	\item a constant feedback $K= K_k$, for $k=0, \dots, N-1$,
	\item a constant covariance $ W = W_k$ for $k = 0, \dots, N-1$,
	\item a constant sensitivity function $G(x,u) = \frac{\partial \psi}{ \partial w}(x,u)$,
	\item the index sets of constraints to be tightened: $\mathcal{I}_{k, \mathrm{tight}}$, which are assumed to be constant for the intermediate shooting nodes $k=1,\dots,N-1$, but might be different for the initial $k=0$ and terminal $k=N$ shooting node,
	\item the backoff factor $\gamma$.
\end{itemize}
The matrices $W_k$ are often constant, which corresponds to an invariant noise distribution.
It is possible to update the $W_k$ matrices outside, based on the current nominal iterates to realize a state-, control-, or time-dependent noise distribution, as e.g. done in~\cite{Lahr2023}.
Moreover, the initial uncertainty matrix $\bar{P}_0$ can be conveniently changed during runtime.
The implementation could be easily modified to allow for different matrices $K_k, W_k$ at each shooting node and respective updates during runtime.

\section{Numerical experiments}
\label{sec:experiments}
This Section presents numerical results comparing the proposed zoRO implementation, labeled \textit{\newzoro}, with previously existing implementations that were tailored to the specific examples in~\cite{Zanelli2021, Gao2023} and are labeled \textit{\oldzoro}.
The experiments have been performed on a Laptop with an Intel i5-8365U CPU, 16~GB of RAM running Ubuntu 22.04.
All code to reproduce the results is open-source.
All experiments have been conducted with \href{https://github.com/acados/acados/releases/tag/v0.2.5}{\acados~v.0.2.5}.

\subsection{Chain benchmark}
The code of the chain benchmark used in~\cite{Zanelli2021} has been adapted to contain the new zoRO implementation \newzoro\footnote{\url{https://github.com/FreyJo/zoro-NMPC-2021}}.
Additionally, the benchmark includes a \textit{nominal} controller, i.e., without uncertainties and back-offs, a \textit{naive} robust version, where the uncertainty is implemented by augmenting the state only exploiting symmetry of the uncertainty matrices, and the zoRO implementation \textit{\oldzoro} proposed in~\cite{Zanelli2021}, which propagates the uncertainties and computes the back-offs in Python.

\begin{figure}[h]
	\includegraphics[width=\linewidth]{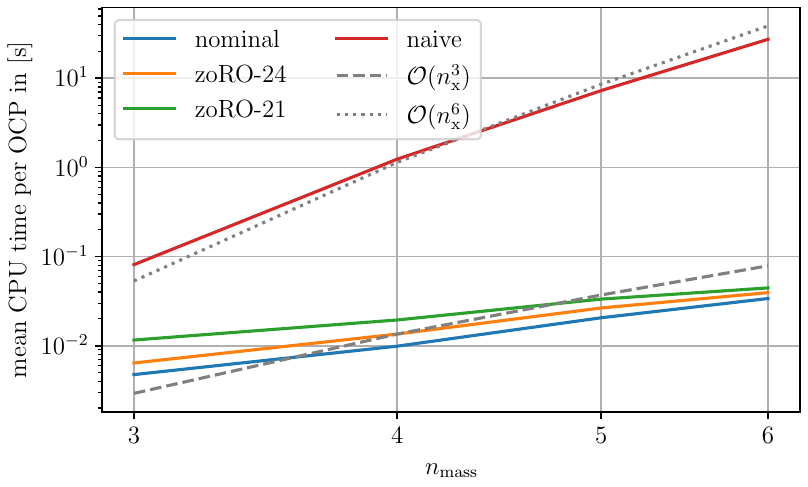}
	\caption{Mean computation time for solving one OCP for different number of masses $n_{\mathrm{mass}}$.
	\label{fig:timings_nmx}
	}
\end{figure}

Figure \ref{fig:timings_nmx} shows a timing comparison of these variants when varying the number of masses.
One can observe that the computation times of both zoRO implementations scale with the same order of complexity.
The \newzoro~implementation is able to bring the computational cost of zoRO much closer to the one of the nominal OCP.
For higher state dimensions, the computation time is dominated by other computations, such as the QP solution.

\subsection{Differential drive robot}
\begin{figure}
	\centering
	\vspace{.2cm}
	\includegraphics[width=\columnwidth]{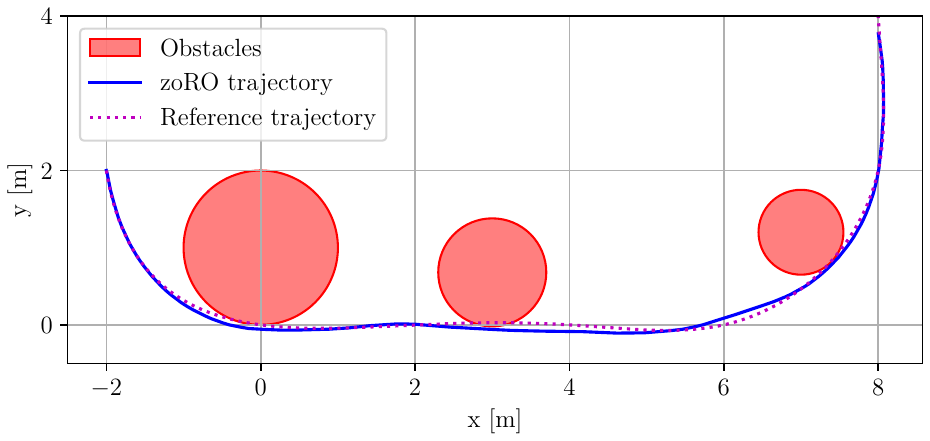}
	\caption{Closed loop trajectory differential drive robot.}
	\label{fig:robo_traj}
\end{figure}

\begin{figure}
	\centering
	\includegraphics[width=\columnwidth]{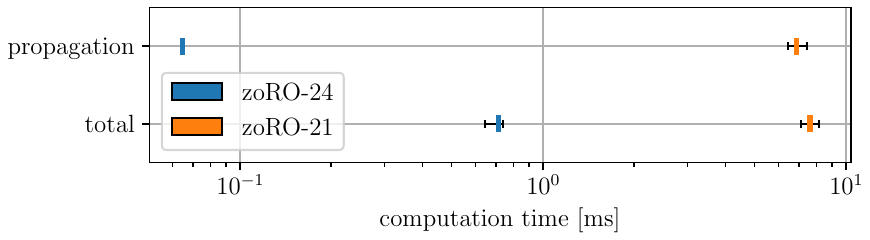}
	\caption{Computation times of zoRO variants. The whiskers indicate the minimum and maximum value.}
	\label{fig:robo_cpu}
\end{figure}

\begin{figure}
	\centering
	\includegraphics[width=\columnwidth]{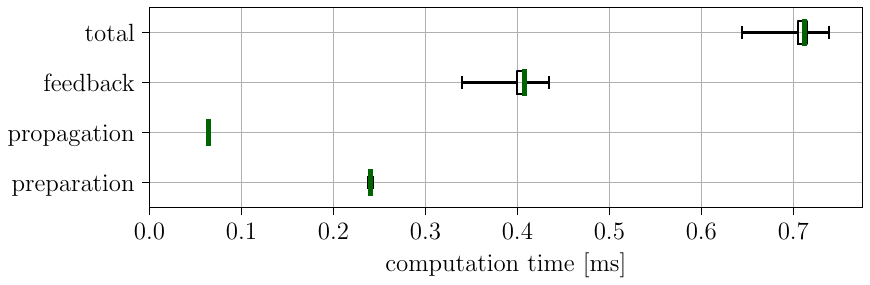}
	\caption{Computation times for the proposed zoRO implementation zoRO-24 on the differential drive robot. The whiskers indicate the minimum and maximum value.}
	\label{fig:robo_cpu_fast_zoro}
\end{figure}
In order to illustrate the effectiveness of our implementation on a practically relevant problem, we regard the system of a differential drive robot, which has been considered in~\cite{Gao2023} and controlled with a zoRO implementation which uses Python for the uncertainty propagation and back-off computation, labeled \oldzoro.
In this system, the state vector is given by $x\!=\!(p_x, p_y, \theta, v, \omega)^\top$, where $p_x, p_y$ parametrize the 2D-position of the center of the robot, $\theta$ is the heading angle, $v$ the forward velocity and $\omega$ the angular velocity.
The controls are the forward acceleration $a$ and angular acceleration $\alpha$.
The ODE is given by
\begin{align}
\dtotal{x}{t} = (v \cos(\theta), v \sin(\theta), \omega, a, \alpha)^\top.
\end{align}
In this OCP, the cost penalizes deviations from a given reference trajectory by linear least-squares.
A multiple shooting discretization with 20 shooting intervals and a horizon of $ 2.0 \mathrm{s} $ is used.
To arrive at discrete dynamics, an implicit Runge-Kutta integrator of order four with Gauss-Legendre Butcher tableau is chosen, which performs three Newton iterations over which the Jacobian is reused.
The constraints consist of bounds on the controls and the states $ v, \omega $, and collision avoidance constraints of the form:
\begin{align}
\norm{(p_x, p_y)\tran - (q_{x,i}, q_{y,i})\tran}_2 \geq r + r_i^{\mathrm{obs}},
\end{align}
where $r$ denotes the radius of the robot, $r_i^{\mathrm{obs}}$ the radius of the obstacle $i$ at position $(q_{x,i}, q_{y,i})$ for $i = 1,\dots, n\ind{obs}$.

The index set of constraints to be tightened $\mathcal{I}_{k, \mathrm{tight}} $ is set to contain the upper bounds on $v$ and $\omega$, the lower bound on $v$ and the collision avoidance constraints.
The terminal constraint consists of very tight bounds on the velocities, corresponding to the nominal trajectory of the robot not moving at the end of the horizon.
Those constraints are not tightened to avoid infeasibility.

The robot is simulated in closed-loop by adding process noise to the output of an RK4 integrator at every simulation instance.
The process noise is sampled from a multivariate normal distribution with zero mean and covariance $W = \diag{2\cdot 10^{-6}, 2\cdot 10^{-6}, 4\cdot 10^{-6}, 1.5\cdot 10^{-3}, 7\cdot 10^{-3}}$.
The initial uncertainty matrix is $\bar{P}_0 = W$, and $\gamma = 3.0$.
The QP subproblems are solved using full condensing and \texttt{DAQP}~\cite{Arnstrom2022}.
The code to reproduce the results discussed next is publicly available\footnote{\url{https://github.com/acados/acados/tree/v0.2.5/examples/acados_python/zoRO_example/diff_drive}}.

Figure~\ref{fig:robo_traj} visualizes the closed-loop trajectory and a reference trajectory of the robot in a scenario with three obstacles.
We observe that the robot is able to avoid the obstacles, even in the described scenario with process noise.
Figure~\ref{fig:robo_cpu} visualizes the computation times of both zoRO implementations.
All timings are obtained by running the exact same simulation 50 times with the same noise realization and taking the minimum of each execution to remove artifacts.
It can be seen that the total computation time of \newzoro~is roughly ten times lower compared to \oldzoro.
The share of computation time for the propagation and back-off computation is $\approx90\%$ for the \oldzoro~implementation and only $\approx10\%$ for \newzoro.

Figure~\ref{fig:robo_cpu_fast_zoro} shows how the computation times of different algorithmic components of \newzoro~on this example.
The total computation time consists of the times corresponding to the \texttt{acados} \textit{preparation} and \textit{feedback} phases, and the uncertainty propagation and backoff computation, labeled as \textit{propagation}.
The preparation step contains the linearization operations of the OCP, such as nonlinear constraints and simulation of the dynamics with sensitivities.
The computation time of the feedback phase is dominated by the \textit{QP} solution.
The propagation and backoff computation can be carried out either as part of the preparation or the feedback phase.
The initial uncertainty matrix $\bar{P}_0$ is constant in this example.
However, if it was obtained from a state estimator, it might be worth to use the new initial uncertainty of the initial state, since the corresponding computations can be carried out quickly relative to the QP solution.

\section{Conclusion \& Outlook}
\label{sec:conclusion}
This paper presented an efficient implementation of the zoRO algorithm for real-time application with an SQP-type solver in \texttt{acados} detailing its stochastic and robust interpretation.
We believe that this work allows the realization of embedded, tube-based MPC in a wider range of real-world control tasks.

The effectiveness and flexibility of the proposed implementation has been demonstrated on two examples from previous studies.
Future work includes the extension of the presented implementation to optionally perform adjoint corrections as proposed in~\cite{Feng2020} and to optimize over linear feedback matrices~\cite{Messerer2021}.

\bibliographystyle{ieeetr}
\bibliography{syscop,amon}
\addcontentsline{toc}{chapter}{Bibliography}

\end{document}